\documentclass{article}

\begin{document}

\title{Conditional Probability formula as a consequence of the Insufficient Reason Principle}

\author{Alexander Dukhovny\\
\small{Department of Mathematics, San Francisco State University} \\
\small{San Francisco, CA 94132, USA} \\
\small{dukhovny [at] sfsu.edu}\\ }

\date{\today}          

\maketitle

\begin{abstract}

The standard conditional probability definition formula is derived as a consequence of the Insufficient Reason Principle expressed as the Maximum Relative Divergence Principle for grading (order-comonotonic) functions on a totally ordered set.

\end{abstract}

\section{Introduction}
\label{sec:Intro}

The Conditional Probability (CP) formula is a pillar of the mathematical approach to many  applications, too numerous to list here (AI, in particular). However, the formula itself has not yet been derived (as far as we know) from a fundamental principle in the way that the "equal likelihood of elementary outcomes" assumption follows from the Insufficient Reason Principle (IRP) in the absence of any information to the contrary.

In numerous probability theory problems IRP has been expressed as the Maximum Shannon Entropy Principle (usually abbreviated as MEP) for choosing the IRP-suggested probability distribution (and/or its missing parameters) as the one that maximizes, under the application-specific constraints, the Shannon Entropy functional (see, e.g., \cite{Shannon}, \cite{Jaynes}). That approach has proved extremely effective in such a large number of cases that references are, again, too many to quote.

The original Shannon Entropy formula and MEP have been generalized in many ways: Kullback-Leibler (K-L) Relative Entropy and Divergence, Partition Entropy, Kolmogorov-Sinai Entropy, Topological Entropy, 
Entropy of General Measures and many others - see, e.g., the review in \cite{Review}. 

In this paper we use the approach started in \cite{Dukh} and continued in \cite{AOME} to generalize Shannon Entropy and Kullback-Leibler entropy concepts to that of the Relative Divergence (RD) of one Grading (order-comonotonic) Function from another on a total order. (The term Relative Divergence comes from Kullback-Leibler Divergence of probability distributions - see \cite{K-L}.) 

Accordingly, Maximum Relative Divergence Principle is used in place of MEP as the mathematical expression of IRP as a method to obtain the IRP-suggested grading function $F$ on a totally ordered set $W$ where a prior (in the Bayesian sense) "null" grading function $G$ is available.

\section{Basic Definitions and Properties}\label{Basics}

The initial RD setup (see \cite{AOME}) begins as follows: let $W = \{ w_k, \quad k \in Z$\} be a chain totally ordered by a relation $\prec$. A real-valued function $F$ on $W$ is said to be a grading function (GF) on $W$ if it is order-comonotonic, that is,

$w \prec v \iff F(w) < F(v)$ for all $w, v \in W $.

(For example, the "indexing" function $I: I(w_k) = k$ is a "natural" GF on $W$.)

For grading functions $F(w)$ and $G(w)$ defined on $W$, the relative divergence of $F$ from $G$ on $W$ is defined as

\begin{equation} \label{RDcountable}
\mathcal{D}(F \Vert G) \vert_W = -\ \sum_{k \in Z}
\ln  \left( \frac{f_k}{g_k}\right) f_k, 
\end{equation}

\noindent where

$f_k=\Delta_k F = F(w_k) - F(w_{k-1})$, $ g_k=\Delta_k G = G(w_k) - G(w_{k-1}), \quad  k \in Z$

are the increments of, respectively, $F$ and $G$ along the chain $W$. (Absolute convergence of the series must be assumed where $W$ is infinite.)

In the special case where $F$ is a probability cumulative distribution function on $W$ and $G = I$, Equation (\ref{RDcountable}) reduces to Shannon Entropy of that probability distribution:

\begin{equation} \label{RD-Shannon}
\mathcal{D}(F \Vert I) \vert_W = -\ \sum_{k \in Z} f_k \ln {f_k}   , 
\end{equation}

Accordingly, the Maximum Relative Divergence Principle (MRDP) for chains is introduced 
as a generalization of the Maximum Entropy Principle as follows:

MRDP: Suppose a "null" grading function $G$ is defined on a chain $W$. Among all application-admissible grading functions on $W$ with the same value range, $F$ is said to be IRP-suggested (or "least-presuming") if its Relative Divergence from the given $G$ is the highest possible.

\section{The CP formula proof}

Using MRDP as an expression of IRP, here is a proof of the classic Conditional Probability formula:

$P(B|A) = \frac{P(A \cap B)}{P(A)}$.

Proof. Let $W$ be a chain of events $W = \{ \emptyset , A\cap B, A \} $ in a random experiment ordered by inclusion. Define function $G$ on $W$ with values, respectively, $\{ 0, p_1 = P(A \cap B), p_2 = P(A)$\} and function $F$ with values $0, x = P(B|A), 1$. As defined, both $G$ and $F$ are grading functions on $W$. 

When $P(A \cap B)$ and $P(A)$ are already known, $G$ must be set as the "null" grading function, whereas $F$ is to be specified by the value of $x$. As such, using MRDP as a mathematical expression of IRP, $x$ must be chosen from the interval $(0,1)$ so as to maximize the value of

$q(x) = \mathcal{D}(F \Vert G) \vert_W.$ 

Now,

$q(x) = - x\ln{ \left( \frac{x}{p_1} \right) } - (1-x)\ln{ \left( \frac{1-x}{p_2 - p_1} \right) },$

$q'(x) = \ln{(1-x)} - \ln{x} + \ln{p_1} - \ln{(p_2 - p_1)}, $

$q''(x) = -{[(1-x)x]^{-1}} < 0$,

so $q(x)$ is concave down on $(0, 1)$ and has the only maximum there at $x = \frac {p_1}{p_2}$, proving the Conditional Probability formula.


\begin{thebibliography} {6}

\bibitem{Shannon}
C.E. Shannon, A Mathematical Theory of Communication. Bell System Technical Journal, vol. 27, July, October (1948), 623-656.

\bibitem{Jaynes}
E.T. Jaynes, Information theory and statistical mechanics, Physical Review. 106 (4) (1957), 620 - 630. 

\bibitem{Dukh}
A. Dukhovny, General Entropy of General Measures, International Journal of Uncertainty, Fuzziness and Knowledge-Based Systems,vol. 10(3) (2002), 213 - 225. 

\bibitem{AOME}
A. Dukhovny,  Axiomatic Origins of Mathematical Entropy: Grading Ordered Sets, arXiv:1903.05240 [math.PR], March 2019

\bibitem{Review}
Jose M. Amigo, Samuel G. Balogh, Sergio Hernandez, A Brief Review of Generalized Entropies, Entropy, 20(2018), 813.  

\bibitem{K-L}
S. Kullback, R.A. Leibler, On information and sufficiency, Annals of Mathematical Statistics. 22 (1)(1951), 79-86. 

\end{thebibliography}
\end{document}